\documentclass[12pt]{article}

 \usepackage{a4}
 \usepackage{amsmath}
 \usepackage{amssymb}
\newtheorem{thm}{Theorem}[section]

 \newtheorem{ex}[thm]{Example}
 \newtheorem{lemma}[thm]{Lemma}
 \newtheorem{rem}[thm]{Remark}

 \DeclareMathOperator{\Rea}{Re}
 
\newcommand{\Int}{\int_{-\infty}^\infty}
 \newcommand{\eps}{\varepsilon}
\newcommand{\Ro}{\mathbb R\setminus\{0\}}
 
 \renewcommand{\a}{\alpha}
  \renewcommand{\b}{\beta}
\title{On powers of Stieltjes moment sequences, II}
 \author{Christian Berg}

 \date{\today}

\begin{document}

 \maketitle

 \begin{abstract}
We consider the set of Stieltjes moment sequences, for which every positive
power is again a Stieltjes moment sequence, we  and prove an integral
representation of the logarithm of the moment sequence in analogy to
the L\'evy-Khinchin representation. We use the result to construct
product convolution semigroups with moments of all orders and to calculate
their Mellin transforms. As an application we construct a positive generating
function for the orthonormal Hermite polynomials.
\end{abstract}



\section{Introduction and main results}
\label{sec:intro}
In his fundamental memoir \cite{St} Stieltjes characterized sequences
of the form
\begin{equation}\label{eq:St}
s_n=\int_0^\infty x^n\;d\mu(x),\;n=0,1,\ldots,
\end{equation}
where $\mu$ is a non-negative measure on $[0,\infty[$,
by certain quadratic forms being non-negative. These sequences are now called
Stieltjes moment sequences. They are called normalized if $s_0=1$. A
 Stieltjes  moment sequence is called {\it S-determinate}, if there is only
  one measure $\mu$ on $[0,\infty[$ such that (\ref{eq:St}) holds;
   ot\-herwise it is
called {\it S-indeterminate}. It is to be noticed that in the S-indeterminate
   case there are also solutions $\mu$ to (\ref{eq:St}), which are not
   supported by $[0,\infty[$, i.e. solutions to the corresponding Hamburger
   moment problem.

A Stieltjes moment sequence is either non-vanishing (i.e. $s_n>0$
    for all $n$) or of the form
  $s_n=c\delta_{0n}$ with $c\geq 0$,
 where $(\delta_{0n})$ is the sequence $(1,0,0,\ldots)$. The latter
 corresponds to the  Dirac measure
  $\delta_0$ with mass $1$ concentrated at $0$.

In this paper we shall characterize the set $\mathcal I$ of normalized
Stieltjes moment
sequences $(s_n)$ with the property that $(s_n^c)$ is a Stieltjes moment
sequence for each $c>0$. The result is given in
 Theorem \ref{thm:harmonic}, from which we extract the following:

 {\it A Stieltjes moment sequence $(s_n)$ belongs to $\mathcal I$
if and only if there exist $\eps\in[0,1]$ and an infinitely divisible probability
$\omega$ on $\mathbb R$  such that
\begin{equation}\label{eq:I}
s_n=(1-\eps)\delta_{0n} + \eps\int_{-\infty}^{\infty} e^{-ny}\,d\omega(y).
\end{equation}
}
We stress however that Theorem \ref{thm:harmonic} also contains a kind of
 L\'evy-Khintchine representation of $\log s_n$ in the case $\eps\neq 0$, and
this result is very useful for deciding if a given sequence belongs 
to $\mathcal I$.

During the preparation of this paper our attention was drawn to the
 Ph.d.-thesis \cite{Ty} of Shu-gwei Tyan, which contains a chapter on
 infinitely divisible moment sequences. The set $\mathcal I$ is the set
 of infinitely divisible Stieltjes moment sequences in the sense of Tyan.
 Theorem 4.2 in \cite{Ty} is a representation
  of $\log s_n$ similar to condition (ii) in Theorem \ref{thm:harmonic}.
As far as we know these results of \cite{Ty} have not been published
elsewhere, so we discuss his results in Section 3.

 The motivation for the present  paper can be found in the paper \cite{B:D}
 by Dur\'an and the author and which provides a unification of recent
 work of Bertoin, Carmona, Petit and Yor, see \cite{B:Y},\cite{C:P:Y1},
 \cite{C:P:Y2}. We found a procedure to generate sequences
 $(s_n)\in\mathcal I$. To formulate the motivating result of
 \cite{B:D}
 we need the concept of a Bernstein function.
                                
Let $(\eta_t)_{t>0}$ be a convolution semigroup of sub-probabilities on
$[0,\infty[$ with {\it Laplace exponent} or {\it Bernstein function}
 $f$ given by
$$
\int_0^\infty e^{-sx}\,d\eta_t(x)=e^{-tf(s)},\;\;s>0,
$$
cf. \cite{B:F},\cite{Bt}. We recall that $f$ has the integral representation
\begin{equation}\label{eq:irepBernstein}
f(s)=a+bs+\int_0^\infty (1-e^{-sx})\,d\nu(x),
\end{equation}
where $a,b\geq 0$ and the L\'evy measure $\nu$ on $]0,\infty[$
satisfies the integrability condition $\int x/(1+x)\,d\nu(x)<\infty$.
Note that $\eta_t([0,\infty[)=\exp(-at)$, so that $(\eta_t)_{t>0}$ consists
of probabilities if and only if $a=0$.

In the following we shall exclude the Bernstein function identically equal to
zero, which corresponds
to the convolution semigroup $\eta_t=\delta_0, t>0$.

 Let $\mathcal B$ denote the set of Bernstein functions which are not
 identically zero. For $f\in\mathcal B$ we note that
  $f'/f$ is completely monotonic as product of the completely monotonic
  functions $f'$ and $1/f$. Therefore  there exists a non-negative
  measure $\kappa$
 on $[0,\infty[$ such that
\begin{equation}\label{eq:f'/f}
 \frac{f'(s)}{f(s)}=\int_0^\infty e^{-sx}\,d\kappa(x).
\end{equation}
It is easy to see that $\kappa(\{0\})=0$ using (\ref{eq:irepBernstein}) and
 $f'(s)\geq \kappa(\{0\})f(s)$. 

\begin{thm}[Berg-Dur\'an \cite{B:D}, Berg \cite{B2}]\label{thm:power}
  Let $\a\geq 0,\b>0$ and let $f\in\mathcal B$ be such that $f(\a)>0$.
 Then the sequence $(s_n)$ defined by 
$$
s_0=1, s_n=f(\a)f(\a+\b)\cdot\ldots\cdot f(\a+(n-1)\b),\quad n\geq 1
$$
belongs to $\mathcal I$.
Furthermore $(s_n^c)$ is S-determinate for $c\leq 2$.
\end{thm}

In most applications of the theorem we put $\a=\b=1$ or $\a=0,\b=1$,
the latter provided $f(0)>0$.
The moment sequence $(s_n^c)$ of Theorem \ref{thm:power} can be
S-indeterminate for $c>2$. This is shown in \cite{B2} for the moment sequences
 \begin{equation}\label{eq:c>2}
 s_n^c=(n!)^c\quad\textrm{and}\quad s_n^c=(n+1)^{c(n+1)}
 \end{equation}
 derived from the Bernstein functions $f(s)=s$ and $f(s)=s(1+1/s)^{s+1}$.
 For the Bernstein function $f(s)=s/(s+1)$ the moment sequence
 $s_n^c=(n+1)^{-c}$ is a Hausdorff moment sequence since
$$
\frac{1}{(n+1)^c}=\frac{1}{\Gamma(c)}\int_0^1 x^n(\log(1/x))^{c-1}\,dx,
$$
and in particular it is S-determinate for all $c>0$. In Section 2 we
give a new proof of Theorem \ref{thm:power}.

In Section 4 we use the Stieltjes moment sequence $(\sqrt{n!})$ to prove
non-negati\-vi\-ty of a generating function for the orthonormal Hermite
 polynomials.

The sequence $(a)_n:=a(a+1)\cdot\ldots\cdot(a+n-1),a>0$ belongs to
$\mathcal I$ and is a one parameter extension of $n!$. For $0<a<b$ also
$(a)_n/(b)_n$ belongs to  $\mathcal I$. These examples are studied in Section
5. Finally, in Section 6 we study a $q$-extension
 $(a;q)_n/(b;q)_n\in\mathcal I$ for $0<q<1,0\leq b<a<1$.
In Section 7 we give some complementary examples.

Any normalized Stieltjes moment sequence $(s_n)$ has the form
$s_n=(1-\eps)\delta_{0n}+\eps t_n$, where $\eps\in[0,1]$ and $(t_n)$ is a normalized
Stieltjes moment sequence satisfying $t_n>0$.

Although the moment sequence $(s_n^c)$ of Theorem \ref{thm:power} can be
S-indeterminate for $c>2$, there is a \lq\lq canonical\rq\rq\, solution
$\rho_c$ to the moment problem defined by \lq\lq infinite divisibility\rq\rq.

The notion of an infinitely divisible probability measure has been studied
for arbitrary locally compact groups, cf. \cite{He}. 

We need the
product convolution $\mu\diamond\nu$ of two measures $\mu$ and $\nu$ on
 $[0,\infty[$: It is
defined as the image of the product measure $\mu\otimes\nu$ under
 the product mapping $(s,t)\mapsto st$. For measures concentrated on
 $]0,\infty[$ it is the convolution with respect to the multiplicative
  group structure on  the interval.
 It is clear that the $n$'th moment of the product convolution is the product
  of the $n$'th moments of the factors.

In accordance with the general definition we say that a probability $\rho$
on $]0,\infty[$ is infinitely divisible on the multiplicative group of
positive real numbers, if
 it has $p$'th product convolution roots for any natural number $p$,
 i.e. if there exists a probability $\tau(p)$ on $]0,\infty[$ such that
 $(\tau(p))^{\diamond p}=\rho$. This condition implies the existence of a
 unique family $(\rho_c)_{c>0}$ of probabilities on
 $]0,\infty[$ such that $\rho_c\diamond\rho_d=\rho_{c+d}$, $\rho_1=\rho$
and $c\mapsto \rho_c$ is weakly continuous. (These conditions imply
that
$\lim_{c\to 0}\rho_c=\delta_1$ weakly.) We call such a family a {\it product
convolution semigroup}. It is a (continuous) convolution semigroup in the
 abstract sense
of \cite{B:F} or \cite{He}. A $p$'th root $\tau(p)$ is unique and one defines
$\rho_{1/p}=\tau(p),\rho_{m/p}=(\tau(p))^{\diamond m}, m=1,2,\ldots.$ Finally
$\rho_c$ is defined by continuity when $c>0$ is irrational.

The family of image measures $(\log(\rho_c))$ under the log-function is a
convolution semigroup of infinitely divisible measures in the ordinary  sense
on the real line considered as an additive group.

The following result generalizes Theorem 1.8 in \cite{B2}, which
treats the special case $\a=\b=1$. In addition we express the Mellin
 transform of the product  
convolution semigroup $(\rho_c)$ in terms of the measure $\kappa$ from
(\ref{eq:f'/f}).

\begin{thm}\label{thm:infdiv} Let $\a\geq 0,\b>0$ and let
  $f\in\mathcal B$ be such that $f(\a)>0$.
The uniquely determined probability measure $\rho$ with moments
$$
s_n=f(\a)f(\a+\b)\cdot\ldots\cdot f(\a+(n-1)\b),\quad n\geq 1
$$
is concentrated
 on $]0,\infty[$ and is
infinitely divisible with respect to
the product convolution.  The unique product convolution semigroup
 $(\rho_c)_{c>0}$  with $\rho_1=\rho$ has the moments
\begin{equation}\label{eq:pcsmoments}
 \int_0^\infty x^n\,d\rho_c(x)=\left(f(\a)f(\a+\b)\cdot\ldots\cdot f(\a+(n-1)\b)\right)^c,\quad c>0, n=1,2,\ldots.
\end{equation}
The Mellin transform of $\rho_c$ is given by
\begin{equation}\label{eq:Mellin}
\int_0^\infty t^z\,d\rho_c(t)=e^{-c\psi(z)},\;\Rea z\geq 0,
\end{equation}
where
\begin{equation}\label{eq:psi}
  \psi(z)=-z\log
  f(\a)+\int_0^\infty\left((1-e^{-z{\b}x})-z(1-e^{-{\b}x})
\right)
  \frac{e^{-{\a}x}}{x(1-e^{-{\b}x})}\,d\kappa(x),
\end{equation}
and $\kappa$ is given by \textrm{(\ref{eq:f'/f})}.
\end{thm}

Proof of the theorem is given in Section 2.

In connection with questions of determinacy the following result is useful.

\begin{lemma}\label{thm:prod} Assume that a Stieltjes moment sequence
 $(u_n)$ is the product
$u_n=s_nt_n$ of two Stieltjes moment sequences $(s_n),(t_n)$. If $t_n>0$ for
 all $n$ and $(s_n)$ is S-indeterminate, then also $(u_n)$ is
 S-indeterminate.
 \end{lemma}

This is proved in Lemma 2.2 and Remark 2.3 in \cite{B:D}.
It follows that if $(s_n)\in\mathcal I$ and $(s_n^c)$ is
S-indeterminate for $c=c_0$, then it is S-indeterminate for any
$c>c_0$. Therefore one of the following three cases occur
\begin{itemize}
\item $(s_n^c)$ is S-determinate for all $c>0$.
\item There exists $c_0, 0<c_0<\infty$ such that $(s_n^c)$ is
  S-determinate  for $0<c<c_0$ and S-indeterminate for $c>c_0$.
\item $(s_n^c)$ is S-indeterminate for all $c>0$.
\end{itemize}
We have already mentioned examples of the first two cases, and the
third case occurs in Remark \ref{thm:BDspecial}.

The question of characterizing the set of normalized Stieltjes moment
 sequences $(s_n)$ with the property that $(s_n^c)$ is a Stieltjes moment
  sequence for each $c>0$ is essentially answered in the monograph
  \cite{B:C:R}. (This was written without knowledge about \cite{Ty}.)
 In fact, $\delta_{0n}$ has clearly this property, so let
  us restrict the attention to the class of non-vanishing normalized
  Stieltjes moment sequences $(s_n)$ for which we can apply the general
   theory of infinitely divisible positive definite kernels, see
   \cite[Proposition 3.2.7]{B:C:R}. Combining this result with Theorem 6.2.6
   in the same monograph we can formulate the solution in the
   following way, where (iii) and (iv) are new:

   \begin{thm}\label{thm:harmonic} For a  sequence
$s_n>0$ the following conditions are equivalent:
\begin{enumerate}
\item[(i)] $s_n^c$ is a normalized Stieltjes moment sequence for each
  $c>0$, i.e. $(s_n)\in\mathcal I$.
\item[(ii)] There exist $a\in\mathbb R, b\geq 0$ and a positive Radon measure
 $\sigma$ on $[0,\infty[\setminus\{1\}$ satisfying
$$
\int_0^\infty (1-x)^2\,d\sigma(x)<\infty,\quad \int_2^\infty x^n\,d\sigma(x)<
\infty,\; n\geq 3
$$
such that
\begin{equation}\label{eq:LKS}
\log s_n=an+bn^2+\int_0^\infty (x^n-1-n(x-1))\,d\sigma(x),\;n=0,1,\ldots.
\end{equation}
\item[(iii)] There exist $0<\eps\leq 1$ and an infinitely divisible probability
$\omega$ on $\mathbb R$ such that
\begin{equation}\label{eq:II}
s_n=(1-\eps)\delta_{0n} + \eps\int_{-\infty}^{\infty} e^{-ny}\,d\omega(y).
\end{equation}
\item[(iv)] There exist $0<\eps\leq 1$ and a product convolution semigroup
$(\rho_c)_{c>0}$ of probabilities on $]0,\infty[$ such that
\begin{equation}\label{eq:III}
s_n^c=(1-{\eps}^c)\delta_{0n}+{\eps}^c\int_0^\infty x^n\,d\rho_c(x),n\geq 0,\;c>0.
\end{equation}
\end{enumerate}
Assume $(s_n)\in\mathcal I$. If $(s_n^c)$ is S-determinate for some
$c=c_0>0$, then the quantities $a,b,\sigma,\eps,\omega,(\rho_c)_{c>0}$
from (ii)-(iv) are uniquely determined. Furthermore $a=\log s_1,b=0$
and the finite measure $(1-x)^2\,d\sigma(x)$ is S-determinate.
\end{thm}

\begin{rem}\label{thm:proof} {\rm The measure $\sigma$ in condition (ii) can
have infinite mass close to 1.
There is nothing special about the lower limit 2 of the second integral. It
can be any number  $>1$.
The conditions on $\sigma$ can be formulated
that $(1-x)^2\,d\sigma(x)$ has moments of any order.
}
\end{rem}

\begin{rem}\label{thm:proof1} {\rm Concerning condition (iv) notice that
the measures
\begin{equation}\label{eq:tilde}
\tilde{\rho}_c =(1-{\eps}^c)\delta_0+{\eps}^c\rho_c,\;c>0
\end{equation}
satisfy the convolution equation
\begin{equation}\label{eq:semi}
 \tilde{\rho}_c\diamond\tilde{\rho}_d=\tilde{\rho}_{c+d}
\end{equation}
  and (\ref{eq:III}) can be written
\begin{equation}\label{eq:IV}
s_n^c=\int_0^\infty x^n\,d\tilde{\rho}_c(x),\;c>0.
\end{equation}
On the other hand, if we start with a family $(\tilde{\rho}_c)_{c>0}$ of
probabilities on $[0,\infty[$ satisfying (\ref{eq:semi}), and if we define
$h(c)=1-\tilde{\rho}_c(\{0\})=\tilde{\rho}_c(]0,\infty[)$, then
$h(c+d)=h(c)h(d)$ and therefore $h(c)={\eps}^c$ with $\eps=h(1)\in
[0,1]$.
 If $\eps=0$
then $\tilde{\rho}_c=\delta_0$ for all $c>0$, and if $\eps>0$ then
$\rho_c:={\eps}^{-c}(\tilde{\rho}_c|]0,\infty[)$ is a probability on $]0,\infty[$
satisfying $\rho_c\diamond\rho_d=\rho_{c+d}$.
}
\end{rem}
 \begin{rem}\label{thm:BDspecial} {\rm In \cite{B:D}  was introduced a
 transformation $\mathcal T$ from normalized
non-vanishing Hausdorff moment sequences $(a_n)$  to normalized Stieltjes
moment sequences $(s_n)$ by the formula
\begin{equation}\label{eq:BDtrans}
\mathcal T[(a_n)]=(s_n);\quad s_n=\frac{1}{a_1\cdot\ldots\cdot a_n},\;n\geq 1.
\end{equation}
We note the following result:

 {\it If $(a_n)$ is a normalized Hausdorff
 moment sequence in $\mathcal I$, then $\mathcal T[(a_n)]\in\mathcal I$.}

 As an example consider the Hausdorff moment sequence $a_n=q^n$, where
$0<q<1$ is fixed. Clearly $(q^n)\in\mathcal I$ and the corresponding
product convolution semigroup is $(\delta_{q^c})_{c>0}$. The transformed
 sequence $(s_n)=\mathcal T[(q^n)]$ is given by
 $$
 s_n=q^{-\binom{n+1}{2}},
 $$
 which again belongs to $\mathcal I$. The sequence $(s_n^c)$ is
 S-indeterminate for all $c>0$. The family of densities
$$
v_c(x)=\frac{q^{c/8}}{\sqrt{2\pi\log(1/q^c)}}\frac{1}{\sqrt{x}}
\exp\left[-\frac{(\log x)^2}{2\log(1/q^c)}\right],\;x>0
$$
form a product convolution semigroup because
$$
\int_0^\infty x^zv_c(x)\,dx=q^{-cz(z+1)/2},\quad z\in\mathbb C.
$$
In particular
 $$
\int_0^\infty x^nv_c(x)\,dx=q^{-c\binom{n+1}{2}}.
 $$
 Each of the measures $v_c(x)\,dx$ is infinitely divisible for the additive
 structure as well as for the multiplicative structure. The additive
 infinite divisibility is deeper than the multiplicative and was first
 proved by Thorin, cf. \cite{Th}.
 }
 \end{rem}

\section{Proofs}
We start by proving Theorem \ref{thm:harmonic} and will deduce Theorem
\ref{thm:power} and \ref{thm:infdiv} from this result.

{\it Proof of Theorem \ref{thm:harmonic}:}
The proof of \lq\lq(i)$\Rightarrow$(ii)\rq\rq is a modification of the proof
 of Theorem 6.2.6 in \cite{B:C:R}: For each $c>0$ we choose a probability
  measure $\tilde{\rho}_c$ on
$[0,\infty[$ such that for $n\geq 0$
 $$
 s_n^c=\int_0^\infty x^n\,d\tilde{\rho}_c(x),
 $$
 hence
 $$
 \int_0^\infty (x^n-1-n(x-1))\,d\tilde{\rho}_c(x)=s_n^c-1-n(s_1^c-1).
 $$
(Because of the possibility of S-indeterminacy we cannot claim the
convolution equation $\tilde{\rho_c}\diamond\tilde{\rho_d}=\tilde{\rho_{c+d}}$.)
 If $\mu$ denotes a vague accumulation point for
  $(1/c)(x-1)^2\,d\tilde{\rho}_c(x)$
 as $c\to 0$, we obtain the representation
 $$
 \log s_n - n\log s_1=\int_0^\infty \frac{x^n-1-n(x-1)}{(1-x)^2}d\mu(x),
 $$
 which gives (\ref{eq:LKS}) by taking out the mass of $\mu$ at $x=1$ and
 defining $\sigma=(x-1)^{-2}d\,\mu(x)$ on $[0,\infty[\setminus\{1\}$.
For details see \cite{B:C:R}.

\medskip
\lq\lq(ii)$\Rightarrow$(iii)\rq\rq\; Define $m=\sigma(\{0\})\geq 0$  and let
$\lambda$ be the image measure on $\Ro$ of
 $\sigma-m\delta_0$ under $-\log x$. We get
 $$
 \int_{[-1,1]\setminus\{0\}} y^2\,d\lambda(y)=\int_{[1/e,e]\setminus\{1\}}
 (1-x)^2\left(\frac{-\log x}{1-x}\right)^2\,d\sigma(x)<\infty,
 $$
and for $n\geq 0$
\begin{equation}\label{eq:Levy}
\int_{\mathbb R\setminus]-1,1[} e^{-ny}\,d\lambda(y)=\int_{]0,\infty[
\setminus]1/e,e[}
x^n\,d\sigma(x)<\infty.
\end{equation}
This shows that $\lambda$ is a L\'evy measure, which allows us to define
a negative definite function
$$
\psi(x)=i\tilde{a}x+bx^2+\int_{\Ro}\left(1-e^{-ixy}-\frac{ixy}{1+y^2}\right)
\,d\lambda(y),
$$
where
$$
\tilde{a}:=\int_{\Ro}\left(\frac{y}{1+y^2}+e^{-y}-1\right)\,d\lambda(y) - a.
$$
Let $(\tau_c)_{c>0}$ be the convolution semigroup on $\mathbb R$ with
$$
\Int e^{-ixy}\,d\tau_c(y)=e^{-c\psi(x)}, x\in\mathbb R.
$$
Because of (\ref{eq:Levy}) we see that $\psi$ and then also $e^{-c\psi}$
 has a holomorphic extension to the lower halfplane.
By a classical result (going back to Landau for Dirichlet series), see
\cite[p.58]{W}, this implies
$$
\Int e^{-ny}\,d\tau_c(y)<\infty, n=0,1,\ldots.
$$
For $z=x+is,s\leq 0$ the holomorphic extension of $\psi$ is given by
$$
\psi(z)=i\tilde{a}z+bz^2+\int_{\Ro}\left(1-e^{-izy}-\frac{izy}{1+y^2}\right)
\,d\lambda(y),
$$
and we have
$$
  \Int e^{-izy}\,d\tau_c(y)=e^{-c\psi(z)}.
$$
In particular we get
\begin{eqnarray*}
-\psi(-in)&=&-\tilde{a}n+bn^2+\int_{\Ro}\left(e^{-ny}-1+\frac{ny}{1+y^2}
\right)\,d\lambda(y)\\
&=&-\tilde{a}n+bn^2+\int_{\Ro}\left(e^{-ny}-1-n(e^{-y}-1)\right)
\,d\lambda(y)\\
&&\quad + n\int_{\Ro}\left(\frac{y}{1+y^2}+ e^{-y}-1\right)\,d\lambda(y)\\
&=& an+bn^2+\int_{]0,\infty[\setminus\{1\}}\left(x^n-1-n(x-1)\right)
\,d\sigma(x),
\end{eqnarray*}
and therefore
\begin{equation}\label{eq:sp}
\log s_n=(n-1)m-\psi(-in) \mathrm{\;for\;} n\geq 1,
\end{equation}
while $\log s_0=\psi(0)=0$.

The measure $\omega=\delta_{-m}*\tau_1$ is infinitely divisible on $\mathbb R$
and we find for $n\geq 1$
$$
s_n=e^{-m}e^{nm-\psi(-in)}=e^{-m}\Int e^{-ny}\,d\omega(y),
$$
so (\ref{eq:II}) holds with $\eps=e^{-m}$.

\medskip
\lq\lq(iii)$\Rightarrow$(iv)\rq\rq\; Suppose (\ref{eq:II}) holds and let
 $(\omega_c)_{c>0}$ be the unique convolution semigroup on $\mathbb R$ such that
 $\omega_1=\omega$. Let $(\rho_c)_{c>0}$ be the product convolution semigroup
 on $]0,\infty[$ such that $\rho_c$ is the image of $\omega_c$ under $e^{-y}$.
 Then (\ref{eq:III}) holds for $c=1,n\geq 0$ and for $c>0$ when $n=0$. For
 $n\geq 1$ we shall prove that
 $$
 s_n^c={\eps}^c\int_0^\infty x^n\,d\rho_c(x),\quad c>0,
 $$
 but this follows from (\ref{eq:II}) first for $c$ rational and then for all
 $c>0$ by continuity.

\medskip
\lq\lq(iv)$\Rightarrow$(i)\rq\rq\; is clear since $(s_n^c)$ is the Stieltjes
 moment  sequence of $\tilde{\rho}_c$ given by (\ref{eq:tilde}).

Assume now $(s_n)\in\mathcal I$. We get $\log s_1=a+b$. If $b>0$ then
$(s_n^c)$ is S-indeterminate for all $c>0$ by Lemma \ref{thm:prod}
because the moment sequence $(\exp(cn^2))$ is S-indeterminate for all
$c>0$ by Remark \ref{thm:BDspecial}.

If $(1-x)^2\,d\sigma(x)$ is S-indeterminate there exist infinitely
many measures $\tau$ on $[0,\infty[$ with $\tau(\{1\})=0$ and such
that
$$
  \int_0^\infty x^n(1-x)^2\,d\sigma(x)= \int_0^\infty
  x^n\,d\tau(x),\quad n\geq 0.
$$
For any of these measures $\tau$ we have
$$
\log s_n=an+bn^2+\int_0^\infty\frac{x^n-1-n(x-1)}{(1-x)^2}\,d\tau(x),
$$
because the integrand is a polynomial. Therefore $(s_n^c)$ has the
S-indeterminate factor
$$
\exp\left(c\int_0^\infty\frac{x^n-1-n(x-1)}{(1-x)^2}\,d\tau(x)\right)
$$
and is itself S-indeterminate for all $c>0$.

We conclude that if $(s_n^c)$ is S-determinate for $0<c<c_0$, then
$b=0$ and $(1-x)^2\,d\sigma(x)$ is S-determinate. Then $a=\log s_1$ and
$\sigma$ is uniquely determined on $[0,\infty[\setminus\{1\}$.
Furthermore, if $\eps,(\rho_c)_{c>0}$ satisfy (\ref{eq:III}) then
$$
   s_n^c=\int_0^\infty x^n\,d\tilde{\rho_c}(x),\quad c>0
$$
with the notation of Remark \ref{thm:proof1}, and we get that
$\tilde{\rho_c}$ is uniquely determined for $0<c<c_0$. This determines
$\eps$ and $\rho_c$ for $0<c<c_0$, but then $\rho_c$ is unique for any
$c>0$ by the convolution equation.

We see that $\eps,\omega$ are uniquely determined by (\ref{eq:II})
since (iii) implies (iv).
$\quad\square$

\bigskip
{\it Proof of Theorem \ref{thm:power} and \ref{thm:infdiv}:}

 To verify directly that the sequence
$$
s_n=f(\a)f(\a+\b)\cdot\ldots\cdot f(\a+(n-1)\b)
$$
 of the form considered in
Theorem \ref{thm:power}
  satisfies (\ref{eq:LKS}), we integrate
formula (\ref{eq:f'/f}) from $\a$ to $s$ and get
 $$
 \log f(s)=\log f(\a) + \int_0^\infty (e^{-{\a}x}-e^{-sx})\frac{d\kappa(x)}{x}.
 $$
Applying this formula we find
 \begin{eqnarray}\label{eq:logsn}
 \log s_n &=& \sum_{k=0}^{n-1} \log f(\a+k\b)\nonumber\\
  &=&n\log f(\a) + \int_0^\infty\left(
 n(1-e^{-{\b}x})-(1-e^{-n{\b}x})\right)\frac{e^{-{\a}x}d\,\kappa(x)}
{x(1-e^{-{\b}x})}\\
 &=&n\log f(\a) + \int_0^1 \left(x^n-1-n(x-1)\right)d\,\sigma(x),\nonumber
 \end{eqnarray}
 where $\sigma$ is the image measure  of
 $$
 \frac{e^{-{\a}x}d\,\kappa(x)}{x(1-e^{-{\b}x})}
 $$
 under $e^{-{\b}x}$. Note that $\sigma$ is concentrated on $]0,1[$.
This shows that $(s_n)\in\mathcal I$.
It follows from the proof of Theorem \ref{thm:harmonic} that the
constant $\eps$ of (iii) is $\eps=1$, so (\ref{eq:III}) reduces to
(\ref{eq:pcsmoments}). The sequence $(s_n^c)$ is S-determinate for
$c\leq 2$ by Carleman's criterion stating that if
\begin{equation}\label{eq:Carleman}
\sum_{n=0}^\infty \frac{1}{\root{2n}\of{s_n^c}}=\infty,
\end{equation}
 then $(s_n^c)$ is S-determinate, cf. \cite{Ak},\cite{S:T}.
To see that this condition is satisfied we note that
$f(s)\leq (f(\b)/\b)s$ for $s\geq \b$, and hence
$$
s_n=f(\a)f(\a+\b)\cdot\ldots\cdot f(\a+(n-1)\b)
$$
$$
\leq
f(\a)(\frac{f(\b)}{\b})^{n-1}\prod_{k=1}^{n-1}(\a+k\b)
=f(\a)f(\b)^{n-1}(1+\frac{\a}{\b})_{n-1}.
$$
It follows from Stirling's formula that (\ref{eq:Carleman}) holds
for $c\leq 2$.

We claim that
\begin{equation}\label{eq:kappa}
\int_1^\infty\frac{e^{-{\a}x}}{x}\,d\kappa(x)<\infty.
\end{equation}
This is clear if $\a>0$, but if $\a=0$ we shall prove
$$
\int_1^\infty\frac{d\kappa(x)}{x}<\infty.
$$
For $\a=0$ we assume that $f(0)=a>0$ and therefore the potential
kernel
$$
p=\int_0^\infty \eta_t\,dt
$$
has finite total mass $1/a$. Furthermore we have
$\kappa=p*(b\delta_0+x\,d\nu(x))$ since
$$
f'(s)=b+\int_0^\infty e^{-sx}x\,d\nu(x),
$$
so we can write $\kappa=\kappa_1+\kappa_2$ with
$$
\kappa_1=p*(b\delta_0+x1_{]0,1[}(x)\,d\nu(x)),\quad
\kappa_2=p*(x1_{[1,\infty[}(x)\,d\nu(x)),
$$
and $\kappa_1$ is a finite measure. Finally
$$
\int_1^\infty\frac{d\kappa_2(x)}{x}=\int_1^\infty\left(\int_0^\infty
\frac{y}{x+y}\,dp(x)\right)d\nu(y)\leq\frac{\nu([1,\infty[)}{a} < \infty.
$$

 The function $\psi$ given by (\ref{eq:psi}) is continuous in the
closed half-plane $\Rea z\geq 0$ and holomorphic in $\Rea z>0$ because
of (\ref{eq:kappa}). Note that
$\psi(n)=-\log s_n$ by (\ref{eq:logsn}). We also notice that $\psi(iy)$ is a
continuous negative definite function on the additive group $(\mathbb R,+)$,
cf. \cite{B:F}, because
$$
  1-e^{-iyx}-iy(1-e^{-x})
$$
is a continuous negative definite function of $y$ for each $x\geq 0$. Therefore
there exists a unique product convolution semigroup $(\tau_c)_{c>0}$
of probabilities on
$]0,\infty[$ such that

\begin{equation}\label{eq:Me1}
   \int_0^\infty t^{iy}\,d\tau_c(t)=e^{-c\psi(iy)},\;c>0,y\in\mathbb R.
\end{equation}

By a classical result, see
\cite[p. 58]{W}), the holomorphy of $\psi$ in the right half-plane implies
that $t^z$ is $\tau_c$-integrable for $\Rea z\geq 0$ and

\begin{equation}\label{eq:Me2}
   \int_0^\infty t^z\,d\tau_c(t)=e^{-c\psi(z)},\;c>0,\Rea z\geq 0.
\end{equation}
In particular the $n$'th moment is given by
$$
  \int_0^\infty t^n\,d\tau_c(t)=e^{-c\psi(n)}=e^{c\log s_n}=s_n^c,
$$
so by S-determinacy of $(s_n^c)$ for $c\leq 2$ we get $\rho_c=\tau_c$ for
$c\leq 2$.This is however enough to ensure that $\rho_c=\tau_c$ for all
$c>0$ since $(\rho_c)$ and $(\tau_c)$ are product convolution semigroups.
$\quad\square$

\bigskip

\section{Tyan's thesis revisited}

In \cite{Ty} Tyan defines a normalized Hamburger moment sequence
 $$
s_n=\Int x^n\mu(x),\quad n\geq 0
$$
 to be  {\it infinitely  divisible} if
\begin{enumerate}
\item[(i)] $s_n\geq 0$ for all $n\geq 0$
\item[(ii)] $(s_n^c)$ is a Hamburger  moment
  sequence for all $c>0$.
\end{enumerate}
Since the set of Hamburger moment sequences is closed under limits and
products, we can replace (ii) by the weaker
\begin{enumerate}
\item[(ii')]$\root{k}\of{s_n}$ is a Hamburger moment sequence for all
  $k=0,1,\ldots$.
\end{enumerate}

\begin{lemma}[Tyan]\label{thm:Tyan1} Let $(s_n)$ be an infinitely divisible
Hamburger moment sequence. Then one of the following cases hold:
\begin{itemize}
 \item $s_n>0$ for all $n$.
\item $s_{2n}>0,s_{2n+1}=0$ for all $n$.
\item $s_n=0$ for $n\geq 1$.
\end{itemize}
\end{lemma}

{\it Proof}: The sequence $(u_n)$ defined by
$$
u_n=\lim_{k\to\infty}\root{k}\of{s_n}=\left\{\begin{array}{ll}
      1 &\mbox{if $s_n>0$} \\
      0 &\mbox{if $s_n=0$}
      \end{array}
      \right.
$$
is a Hamburger moment sequence, and since it is bounded by 1 we have
$$
u_n=\int_{-1}^1 x^n\,d\mu(x)
$$
for some probability $\mu$ on $[-1,1]$.

Either $u_2=1$ and then $\mu=\a\delta_1+(1-\a)\delta_{-1}$ for some
$\a\in[0,1]$, or $u_2=0$ and then $\mu=\delta_0$, which gives the
third case of the Lemma.

In the case $u_2=1$ we have $u_1=2\a-1$, which is either 1 or 0
corresponding to either $\a=1$ or $\a=\frac12$, which gives the two
first cases of the Lemma. $\quad\square$

The symmetric case $s_{2n}>0,s_{2n+1}=0$ is equivalent to studying
infinitely divisible Stieltjes moment sequences, while the third case
is trivial.

Theorem 4.2 of \cite{Ty} can be formulated:
\begin{thm}\label{thm:Tyan2} A Hamburger moment sequence $(s_n)$ such
  that $s_n>0$ for all $n$ is infinitely divisible if and only if
the following representation holds
$$
\log s_n=an+bn^2+\Int(x^n-1-n(x-1))\,d\sigma(x),\quad n\geq 0,
$$
where $a\in\mathbb R,b\geq 0$ and $\sigma$ is a positive measure
on $\mathbb R\setminus\{1\}$ such that $(1-x)^2\,d\sigma(x)$ is a
 measure with moments of any order. Furthermore
$(s_n)$ is a Stieltjes moment sequence if and only if $\sigma$ can  be
chosen supported by $[0,\infty[$.
 \end{thm}

The proof is analogous to the proof of Theorem \ref{thm:harmonic}. 

Tyan also discusses infinitely divisible multidimensional moment
sequences and obtains analogous results. 
\bigskip

\section{An application to Hermite polynomials}

It follows from equation (\ref{eq:c>2}) that
\begin{equation}\label{eq:n!}
\sqrt{n!}=\int_0^\infty u^n\,d\sigma(u)
\end{equation}
for the unique probability $\sigma$ on the half-line satisfying
 $\sigma\diamond\sigma=\exp(-t)1_{]0,\infty[}(t)\,dt$. Even though
 $\sigma$
is not explicitly known, it can be used to prove that a certain
generating function for the Hermite polynomials is non-negative.

Let $H_n,n=0,1,\ldots$ denote the sequence of Hermite polynomials
satisfying the orthogonality relation
$$
\frac{1}{\sqrt{\pi}}\int_{-\infty}^\infty
H_n(x)H_m(x)e^{-x^2}\,dx=2^nn!\delta_{nm}.
$$
The following generating function is well known:
\begin{equation}\label{eq:Hermite}  
\sum_{k=0}^\infty \frac{H_k(x)}{k!}z^k=e^{2xz-z^2},\quad x,z\in\mathbb
C.
\end{equation}

The corresponding orthonormal polynomials are given by
$$
h_n(x)=\frac{H_n(x)}{\sqrt{2^nn!}},
$$
and they satisfy the following inequality of Szasz, cf. \cite{Sz}
\begin{equation}\label{eq:Sz}
|h_n(x)|\leq e^{x^2/2}, x\in\mathbb R, n=0,1,\ldots.
\end{equation}

Let $\mathbb D$ denote the open unit disc in the complex plane.

\begin{thm}\label{thm:Hermite}
The generating function
\begin{equation}\label{eq:gen}
G(t,x)=\sum_{k=0}^\infty h_k(x)t^k
\end{equation}
is continuous for $(t,x)\in \mathbb D\times\mathbb R$ and satisfies
 $G(t,x)>0$ for $-1<t<1,x\in\mathbb R$.
\end{thm}

{\it Proof:} The series for the generating function (\ref{eq:gen})
converges uniformly on compact subsets of $\mathbb D\times\mathbb R$
by the inequality of Szasz (\ref{eq:Sz}), so it is continuous.

By (\ref{eq:n!}) we find  
$$
\sum_{k=0}^n h_k(x)t^k=\int_0^\infty\left(\sum_{k=0}^n
 \frac{H_k(x)}{k!}(\frac{tu}{\sqrt{2}})^k\right)\,d\sigma(u),
$$
which by (\ref{eq:Hermite}) converges to
$$
\int_0^\infty \exp(\sqrt{2}tux-t^2u^2/2)\,d\sigma(u)>0\;
\textrm{for}\; -1<t<1,x\in\mathbb R, 
$$
provided we have dominated convergence. 
This follows however from (\ref{eq:Sz})
because
$$
\int_0^\infty \left|\sum_{k=0}^n
 \frac{H_k(x)}{k!}(\frac{tu}{\sqrt{2}})^k
\right|\,d\sigma(u)\leq e^{x^2/2}\int_0^\infty
\left(\sum_{k=0}^\infty\frac{(|t|u)^k}{\sqrt{k!}}\right)\,d\sigma(u)
$$
$$
=e^{x^2/2}(1-|t|)^{-1}<\infty.
$$

$\quad\square$

\section{The moment sequences $(a)_n^c$ and $((a)_n/(b)_n)^c$}

For each $a>0$ the sequence $(a)_n:=a(a+1)\cdot\ldots\cdot(a+n-1)$ is the
Stieltjes moment sequence  of the $\Gamma$-distribution $\gamma_a$:
$$
(a)_n=\frac{\Gamma(a+n)}{\Gamma(a)}=\int x^n\,d\gamma_a(x)=\frac{1}{\Gamma(a)}
\int_0^\infty x^{a+n-1}e^{-x}\,dx.
$$
For $a=1$ we get the moment sequence $n!$, so the following result generalizes
Theorem 2.5 of \cite{B2}.

\begin{thm}\label{thm:(a)_n^c} The sequence $(a)_n$ belongs to $\mathcal I$ for
each $a>0$. There exists a unique product convolution
semigroup $(\gamma_{a,c})_{c>0}$ such that $\gamma_{a,1}=\gamma_a$. The
moments are given as
$$
\int_0^\infty x^n\,d\gamma_{a,c}(x)=(a)_n^c,\quad c>0,
$$
and
$$
\int_0^\infty
x^z\,d\gamma_{a,c}(x)=\left(\frac{\Gamma(a+z)}{\Gamma(a)}
\right)^c,\quad \Rea z>-a.
$$
The moment sequence $((a)_n^c)$ is S-determinate for $c\leq 2$ and
S-indeterminate for $c>2$.
\end{thm}

{\it Proof:} We apply Theorem \ref{thm:power} and \ref{thm:infdiv} to
the Bernstein function $f(s)=a+s$ and put $\a=0,\b=1$. The formula for
the Mellin transform follows from a classical formula about
$\log\Gamma$, cf. \cite[8.3417]{G:R}.

We shall prove that $(a)_n^c$ is S-indeterminate for $c>2$.
In \cite{B2} it was proved that $(n!)^c$ is S-indeterminate for $c>2$, and
so are all the shifted sequences $((n+k-1)!)^c, k\in\mathbb N$. This implies
that
$$
(k)_n^c=\left(\frac{(n+k-1)!}{(k-1)!}\right)^c
$$
is S-indeterminate for $k\in\mathbb N,c>2$. To see that also $(a)_n^c$ is
S-indeterminate for $a\notin \mathbb N$, we choose an integer $k>a$ and
factorize
$$
(a)_n^c=\left(\frac{(a)_n}{(k)_n}\right)^c (k)_n^c.
$$
By the following theorem the first factor is a non-vanishing Stieltjes moment
sequence, and by Lemma \ref{thm:prod} the product is S-indeterminate.
$\quad\square$

\medskip
For $0<a<b$ we have
\begin{equation}\label{eq:beta}
  \frac{(a)_n}{(b)_n}=\frac{1}{B(a,b-a)}\int_0^1 x^{n+a-1}(1-x)^{b-a-1}\,dx,
  \end{equation}
  where $B$ denotes the Beta-function.

\begin{thm}\label{thm:(a/b)} Let $0<a<b$. Then $((a)_n/(b)_n)$ belongs to
$\mathcal I$ and all powers of the moment sequence are Hausdorff moment
sequences. There exists a unique product convolution semigroup
 $(\b(a,b)_c)_{c>0}$ on $]0,1]$ such that
$$
  \int_0^1 x^z\,d\b(a,b)_c(x)=\left(\frac{\Gamma(a+z)}{\Gamma(a)}/
\frac{\Gamma(b+z)}{\Gamma(b)}\right)^c,\quad \Rea z>-a.
$$
\end{thm}

{\it Proof:} We apply Theorem \ref{thm:power} and \ref{thm:infdiv} to
the Bernstein function $f(s)=(a+s)/(b+s)$ and put $\a=0,\b=1$.

 The Stieltjes moment sequences $(((a)_n/(b)_n)^c))$ are
all bounded and hence Hausdorff moment sequences. The measures
$\gamma_{b,c}\diamond\b(a,b)_c$ and $\gamma_{a,c}$ have the same
moments and are therefore equal for $c\leq 2$ and hence for any $c>0$
by the convolution equations. The Mellin transform of $\b(a,b)_c$
follows from Theorem \ref{thm:(a)_n^c}.

 $\square$

\bigskip
\section{The $q$-extension $((a;q)_n/(b;q)_n)^c$}

In this section we fix $0<q<1$ and consider the $q$-shifted factorials
$$
 (z;q)_n=\prod_{k=0}^{n-1} (1-zq^k), z\in\mathbb C, n=1,2,\ldots,\infty
 $$
 and $(z;q)_0=1$. We refer the reader to \cite{Ga:Ra} for further
 details about $q$-extensions of various functions.

For $0\leq b<a<1$ the sequence $s_n=(a;q)_n/(b;q)_n$ is a Hausdorff moment
sequence for the measure
\begin{equation}\label{eq:q-beta}
\mu(a,b;q)=\frac{(a;q)_\infty}{(b;q)_\infty}\sum_{k=0}^\infty
\frac{(b/a;q)_k}{(q;q)_k}a^k\delta_{q^k},
\end{equation}
which is a probability on $]0,1]$ by the $q$-binomial Theorem, cf.
 \cite{Ga:Ra}.
The calculation of the $n$'th moment follows also from this
 theorem since
 $$
 s_n(\mu(a,b;q))=\frac{(a;q)_\infty}{(b;q)_\infty}\sum_{k=0}^\infty
\frac{(b/a;q)_k}{(q;q)_k}a^kq^{kn}
=\frac{(a;q)_\infty}{(b;q)_\infty}\frac{((b/a)aq^n;q)_\infty}{(aq^n;q)_\infty}
=\frac{(a;q)_n}{(b;q)_n}.
$$

Replacing $a$ by $q^a$ and $b$ by $q^b$ and letting $q\to 1$ we get the
moment sequences $(a)_n/(b)_n$, so the present example can be thought
of as a $q$-extension of the former.  The distribution $\mu(q^a,q^b;q)$ is
 called the $q$-Beta law in Pakes \cite{P1} because of
 its relation to the $q$-Beta function.

\begin{thm}\label{thm:q-betapowers} For $0\leq b<a<1$ the sequence
$s_n=(a;q)_n/(b;q)_n$ belongs to $\mathcal I$.
The measure $\mu(a,b;q)$ is infinitely divisible with respect to the product
 convolution and the corresponding product convolution semigroup
 $(\mu(a,b;q)_c)_{c>0}$ satisfies
\begin{equation}\label{eq:q-betamellin}
\int t^z\,d\mu(a,b;q)_c(t)=\left(\frac{(bq^z;q)_\infty}{(b;q)_\infty}/
\frac{(aq^z;q)_\infty}{(a;q)_\infty}\right)^c,\quad\Rea z>-\frac{\log a}{\log q}.
\end{equation}
In particular
\begin{equation}\label{eq:q-betapowers}                                      
s_n^c=((a;q)_n/(b;q)_n)^c
\end{equation}
is the moment sequence of $\mu(a,b;q)_c$, which is concentrated on
$\{q^k\mid k=0,1,\ldots\}$ for each $c>0$.
\end{thm}

{\it Proof}: It is easy to prove that $(a;q)_n/(b;q)_n$ belongs to $\mathcal I$
using Theorem \ref{thm:power} and \ref{thm:infdiv} applied to
the Bernstein function 
$$
f(s)=\frac{1-aq^s}{1-bq^s}=1-(a-b)\sum_{k=0}^\infty b^kq^{(k+1)s},
$$
 but it will also be a consequence of
the following considerations, which gives information about the
support of $\mu(a,b;q)_c$.

 For a probability $\mu$ on $]0,1]$ let $\tau=-\log(\mu)$ be
the image measure of $\mu$  under $-\log$. It is concentrated on $[0,\infty[$
and
$$
\int_0^1 t^{ix}\,d\mu(t)=\int_0^\infty e^{-itx}\,d\tau(t).
$$
This shows that $\mu$ is infinitely divisible with respect to the product
convolution if and only if $\tau$ is infinitely divisible in the ordinary
sense, and in the affirmative case the negative definite function $\psi$
associated to $\mu$ is related to the Bernstein function $f$ associated to
$\tau$ by $\psi(x)=f(ix), x\in\mathbb R$, cf. \cite[p.69]{B:F}.

We now prove that $\mu(a,b;q)$ is infinitely divisible for
the product convolution. As noticed this is equivalent to proving that the
 measure
$$
\tau(a,b;q):=\frac{(a;q)_\infty}{(b;q)_\infty}\sum_{k=0}^\infty
\frac{(b/a;q)_k}{(q;q)_k}a^k\delta_{k\log(1/q)},
$$
is infinitely divisible in the ordinary sense. To see this we calculate the
Laplace transform of $\tau(a,b;q)$ and get by the  $q$-binomial Theorem

\begin{equation}\label{eq:Laplace}
\int_0^\infty e^{-st}\,d\tau(a,b;q)(t)= \frac{(bq^s;q)_\infty}{(b;q)_\infty}/
\frac{(aq^s;q)_\infty}{(a;q)_\infty},\quad s\geq 0.
\end{equation}

Putting
$$
f_a(s)=\log \frac{(aq^s;q)_\infty}{(a;q)_\infty},
$$
we see that $f_a$ is a bounded Bernstein function  of the form
$$
f_a(s)=-\log(a;q)_\infty-\varphi_a(s),
$$
where
$$
\varphi_a(s)=-\log(aq^s;q)_\infty=\sum_{k=1}^\infty \frac{a^k}{k(1-q^k)}q^{ks}
$$
is completely monotonic as Laplace transform of the finite measure
$$
\nu_a=\sum_{k=1}^\infty\frac{a^k}{k(1-q^k)}\delta_{k\log(1/q)}.
$$
 From (\ref{eq:Laplace}) we get
 $$
 \int_0^\infty
 e^{-st}\,d\tau(a,b;q)(t)=\frac{(a;q)_\infty}{(b;q)_\infty}
e^{\varphi_a(s)-\varphi_b(s)},
 $$
and it follows that $\tau(a,b;q)$ is infinitely divisible and the
corresponding convolution semigroup is given by the infinite series
$$
\tau(a,b;q)_c=\left(\frac{(a;q)_\infty}{(b;q)_\infty}\right)^c\sum_{k=0}^\infty
\frac{c^k(\nu_a-\nu_b)^{*k}}{k!},\quad c>0.
$$

Note that each of these measures are concentrated on
$\{k\log(1/q)\mid k=0,1,\ldots\}$. The associated L\'evy measure is
the finite measure $\nu_a-\nu_b$ concentrated on $\{k\log(1/q)\mid
k=1,2,\ldots\}$. This shows that the image measures
 $$
\mu(a,b;q)_c=\exp(-\tau(a,b;q)_c),\quad c>0
$$
 form a product convolution semigroup
concentrated on $\{q^k\mid k=0,1,\ldots\}$.

The product convolution semigroup $(\mu(a,b;q)_c)_{c>0}$ has the negative
 definite function $f(ix)$, where $f(s)=f_a(s)-f_b(s)$ for
 $\Rea s\geq 0$, hence
$$
\int t^{ix}\,d\mu(a,b;q)_c(t)=\left(\frac{(bq^{ix};q)_\infty}{(b;q)_\infty}/
\frac{(aq^{ix};q)_\infty}{(a;q)_\infty}\right)^c,\quad x\in\mathbb R,
$$
and the equation (\ref{eq:q-betamellin}) follows by holomorphic continuation.
Putting $z=n$ gives (\ref{eq:q-betapowers}).

$\quad\square$

\bigskip
\section{Complements}
\label{sec:com}

\begin{ex}
{\rm  Let $0<a<b$ and consider the Hausdorff moment sequence
$a_n=(a)_n/(b)_n\in\mathcal I.$ By Remark \ref{thm:BDspecial} the
moment sequence $(s_n)=\mathcal T[(a_n)]$ belongs to $\mathcal I$. We find
$$
s_n=\prod_{k=1}^n\frac{(b)_k}{(a)_k}=\prod_{k=0}^{n-1}
\left(\frac{b+k}{a+k}\right)^{n-k}.
$$
}
\end{ex}

\begin{ex}
{\rm
Applying $\mathcal T$ to the Hausdorff moment sequence
$((a;q)_n/(b;q)_n)$ gives the Stieltjes moment sequence
\begin{equation}\label{eq:q-binom}
s_n=\prod_{k=1}^n \frac{(b;q)_k}{(a;q)_k}=
\prod_{k=0}^{n-1}\left(\frac{1-bq^k}{1-aq^k}\right)^{n-k}.
\end{equation}

We shall now give the measure with moments (\ref{eq:q-binom}).

For $0\leq p<1, 0<q<1$ we consider the function of $z$
$$
h_p(z;q)=\prod_{k=0}^\infty\left(\frac{1-pzq^k}{1-zq^k}\right)^k,
$$
which is holomorphic in the unit disk with a power series expansion
\begin{equation}\label{eq:power}
h_p(z;q)=\sum_{k=0}^\infty c_kz^k
\end{equation}
having non-negative coefficients $c_k=c_k(p,q)$. To see this, notice that
$$
\frac{1-pz}{1-z}=1+ \sum_{k=1}^\infty (1-p)z^k.
$$
For $0\leq b<a<1$ and $\gamma>0$ we consider the probability measure with
support in $[0,\gamma]$
$$
\sigma_{a,b,\gamma}=\frac{1}{h_{b/a}(a;q)}\sum_{k=0}^\infty c_ka^k\delta_{{\gamma}q^k},
$$
where the numbers $c_k$ are the (non-negative) coefficients of the power
 series for $h_{b/a}(z;q)$.

The $n$'th moment of $\sigma_{a,b,\gamma}$ is given by
$$
s_n(\sigma_{a,b,\gamma})={\gamma}^n\frac{h_{b/a}(aq^n;q)}{h_{b/a}(a;q)}.
$$

For $\gamma=(b;q)_\infty/(a;q)_\infty$ it is easy to see that
$$
s_n(\sigma_{a,b,\gamma})=\prod_{k=0}^{n-1}
\left(\frac{1-bq^k}{1-aq^k}\right)^{n-k},
$$
which are the moments (\ref{eq:q-binom}).
}
\end{ex}

\author{C. Berg, Department of Mathematics, University of Copenhagen,
  Universitetsparken 5, DK-2100, Denmark; Email: berg@math.ku.dk 

\end{document}